# Much ado about zero

Asis Kumar Chaudhuri
Variable Energy Cycltron Centre
Kolkata-700 064

**Abstract:** A brief historical introduction for the enigmatic number Zero is given. The discussions are for popular consumption.

Zero is such a well known number today that to talk about it may seem "much ado about nothing." However, in the history of mathematics, the concept of zero is rather new. Homo sapiens or humans appeared on the earth approximately 200,000 years ago and only around 3000 BCE, men learned mathematics and only around 628 CE, zero was discovered in India. There is a saying that "necessity is the mother of inventions". Out of necessity only, men learned mathematics. Over the ages, initially nomadic men learned to light the fire, learned to herd animals, learned agriculture and settled around the fertile river valleys, where the essential ingredient for farming, water; was available in plenty. All the ancient civilizations were river bank civilizations e.g. Mesopotamian civilization (part of modern Iraq) between Tigris and Euphrates, Indus valley civilization (now in Pakistan) along the river Indus, Egyptian civilization along the Nile river and Huan He civilization along the river Huan He in China. Men settled with agriculture and husbandry needed to learn mathematics. It was required to count their domesticated animals, to measure the plot of land, to fix taxation, etc. It is true that intuitively, from the very beginning, men could distinguish between, say one horse and two horses, but could not distinguish between say 25 horses and 26 horses. Even now, in several tribal societies, counting beyond three or four is beyond their ability. For them anything more than three or four is many. The following story from George Gamow's famous book, "One Two Three Infinity" amply illustrate the situation,

*Two Hungarian aristocrats decided to play a game in which the one who calls the largest number wins.*

*"Well," said one of them, "you name your number first."*

*After a few minutes of hard mental work the second aristocrat finally named the largest number he could think of.*

*"Three," he said.*

*Now it was the turn of the first one to do the thinking, but a quarter of an hour he finally gave up.*

*"You've won," he agreed.*

Sumer, in the Mesopotamian region is called the cradle of civilization. Sumerians were the first people to introduce plow for farming, introduced wheel for movement. They were also the first people to learn writing. They wrote on clay tablets with wedge shaped implements and their script is called "cunieform". Sumerians were also the first people to learn mathematics and used symbols to represent numbers. They learned to count large numbers. Counting is a human action which requires certain mental developments, which took many ages to come. Indeed, we have been tuned with the idea; otherwise, that the symbol "3" can be used to represent three horses as well as three cows would have been astonishing to us. While we understand 3 horses or 3 cows, but, what exactly 3 represents in mathematics is not easily understood. It required great mathematician and logician like Gottlob Frege (8 November 1848 – 26 July 1925). Giuseppe Peano (1858- 1932), and Bertrand Russell (18 May 1872 – 2 February 1970) to tell us its true meaning, that too only after the concept of set was introduced by George Cantor (March 3 1845 – January 6, 1918) and Richar Dedekind (6 October 1831 – 12 February 1916). Simply speaking, a set is a collection; collection of anything. For example, the shirt, pant, shoes, underwear, you are wearing can form a collection or a set. One can also define an empty set or a null set that contains absolutely nothing. With the help of the set theory, number 3 is defined as the set whose members each have 3 elements. The definition could be generalized, a natural number n is the set whose members each have n elements. Rather an abstract definition!

Over the years, the city state Sumer, expanded to the Babylonian Empire. Babylonians used a positional number system (where symbols positions indicate their value) yet did not use zero. Without textual help they could not distinguish, say between, 312 and 3012. Like the Sumerians, Babylonians also used sexagesimal (base 60) numeral system. In a sexagesimal system, 60 distinct symbols are used to represent the numbers 1,2,3...,60. Truly speaking, theirs were not a pure sexagesimal system as they used fewer symbols and combined them to count up to 60. One may wonder how they came up with the base 60 system. There are some suggestions that the origin of sexagesimal system may be related to early man's calculation of 360 days a year. To quote from Florian Cajori's book, "A History of Mathematics",

*"Cantor offers the following theory: At first the Babylonians reckoned the year at 360 days. This led to the division of the circle into 360 degrees, each degree representing the daily amount of the supposed yearly revolution of the sun around the earth. Now they were, very probably, familiar with the fact that the radius can be applied to its circumference as a chord 6 times, and that each of these chords subtends an arc measuring exactly 60 degrees. Fixing their attention upon these degrees, the division into 60 parts may have suggested itself to them. Thus, when greater precision necessitated a subdivision of the degree, it was partitioned into 60 minutes. In this way the sexagesimal notation may have originated."*

I am not specifically fond of this theory. To me, it presupposes rather an advanced knowledge, that the radius of a circle can be applied to its circumference as a chord 6 times. I will personally favor the alternate suggestion that it is easy to count up to 60 with 5 fingers of one hand and 12 knuckles of the other. There was one more compelling reason. Sumerian's did not have the radix point, the symbol we use to separate the integer part from the fractional part of a number. 60 is a super composite number. It is divided by; 1,2,3,4,5,6,10,12,15,20, 30 and 60. Use of fractions is the minimized in sexagesimal or base 60 system. A Remnant of the base 60 system can still be found in our division of time, e.g. One hour is 60 minutes, one minute is 60 seconds. Sumerians also developed a method to describe very large numbers. While Sumerian and Babylonian mathematics were quite developed, they did not have zero. Over the years, Babylonians developed a method, whereby a vacant space is left in between two numbers to indicate a place without value. While some European historians likened the "vacant space" with modern zero, but truly, it is nothing like zero. The vacant space is always between two numbers, never at the end of a number, the way we designate 30 or 300. Later, the Sumerian knowledge of mathematics passed on to other ancient civilizations e.g. Egypt, India, and China etc. All these civilizations contributed to the growth of mathematics. For example, while Sumerians and Babylonians used sexagesimal numeric system, Egyptian and Indians used the decimal or base 10 numeric system, the system which is prevalent today. Possibly, Egyptian's invented the base 10 system for ease of counting up to 10 using 10 fingers of the two hands.

The Greeks, who contributed vastly in geometry and number theory, adopted a system based on the additive principle and had no need to introduce zero. They thought of numbers geometrically, for example in Euclid's Element, there are fine discussions on the number theory, but numbers are always represented as lengths of a line. Greek's did not have separate symbols for numbers; rather they were drawn from their list of alphabets. Few of the Greek numerals are shown in the figure 1. They had separate symbols for higher numbers, e.g. 100=ρ, 500=φ, etc. The highest number they could write is 10000=μ, called "Myrad".

| alpha | beta | gamma | delta | epsilon | stigma | zeta | eta | theta | iota |
|-------|------|-------|-------|---------|--------|------|-----|-------|------|
| α | β | γ | δ | ε | ς | ξ | η | θ | ι |
| 1 | 2 | 3 | 4 | 5 | 6 | 7 | 8 | 9 | 10 |

*Figure 1: Greek's draw their numerals from the list of alphabets. Few of the Greek numerals are shown here.*

Roman number system (still in limited use) is similar to that of Greeks. Seven symbols;

I = 1, V = 5, X = 10, L = 50, C = 100, D = 500, and M = 1,000,

are used to construct all sort of numbers. In this system a letter placed after another of greater value adds and letter placed before another of greater value subtracts. Thus;

III=3, IX=9 and XI=11, MMXVI=2016.

You can imagine the grief of a Roman mathematician if he has to write a large number. For example, to write 100 thousand, he has to write, M one after another, one hundred times.

Zero was invented in India. Contribution of India, in the development of mathematics, can be best appreciated from the following quote from Laplace[1] (1749- 1827), one of the greatest mathematician from France,

*"It is India that gave us the ingenious method of expressing all numbers by means of ten symbols, each symbol receiving a value of position as well as an absolute value; a profound and important idea which appears so simple to us now that we ignore its true merit. But it's very simplicity and the great ease which it has lent to computations put our arithmetic in the first rank of useful inventions; and we shall appreciate the grandeur of the achievement the more when we remember that it escaped the genius of Archimedes[2] and Apollonius[3], two of the greatest men produced by antiquity"*

Zero was invented only around 628 CE, but Indian mathematics and astronomy reached excellence much before that. Traditionally, Indians used Brahmi numeral system, the first ten numbers are shown in Fig.2. It was a decimal system but did not have a zero. It was also not a place value system and separate symbols were used for 20, 30, 40 etc.

| 1 | 2 | 3 | 4 | 5 | 6 | 7 | 8 | 9 | 10 |
|---|---|---|---|---|---|---|---|---|---|
| − | = | ≡ | + | Һ | Ψ | ? | ﻼ | ? | ∝ |

*Figure 2: Brahmi numeral system.*

---

[1] *Students of mathematics and physics are well acquainted with the French mathematician Pierre-Simon marquis de Laplace (1749- 1827CE). Laplace transform. Laplacian differential operator, widely used in mathematics, is named after him. He pioneered the probability theory. He also developed the nebular hypothesis of the origin of the Solar System and was one of the first scientists to postulate the existence of black holes and the notion of gravitational collapse.*

[2] *Archimedes, (ca. 290 BCE- 212/211 BCE), was the most famous mathematician and inventor in ancient Greek. His achievements were numerous. He discovered, Lever, pulley, screw pump, war machinery etc. However, he is most known for his formulation of the hydrostatic principle (also known as Archimedes' principle).*

[3] *Apollonius of Perga, (ca 240 BCE-ca 190 BCE) was a Greek mathematician and Geometer, known for his works on 'Conic sections. The modern name ellipse, parabola, hyperbola are due to him. He greatly influenced Ptolemy, Copernicus, Newton and others.*

Aryabhata, India's most important post-Vedic mathematician, and astronomer lived during 476-550 CE. He wrote several mathematical and astronomical treatises, most of which are now lost. His most important work, Aryabhatiya had survived up to the modern times and was written when he was only 23 years old. The first chapter of Aryabhatia is Dasagitika. In Dasagitika, Aryabhata, possibly under the influence of the Greeks, introduced a numeral system based on the Sanskrit alphabets.

| क =1 | ख =2 | ग =3 | घ =4 | ङ =5 |
|---|---|---|---|---|
| च =6 | छ =7 | ज =8 | झ =9 | ञ =10 |
| ट =11 | ठ =12 | ड =13 | ढ =14 | ण =15 |
| त =16 | थ =17 | द =18 | ध =19 | न =20 |
| प =21 | फ =22 | ब =23 | भ =24 | म =25 |
| य =30 | र =40 | ल =50 | व =60 | |
| श =70 | ष =80 | स =90 | | |
| ह =100 | | | | |

*Figure 3. Numeric values assigned by Aryabhata to Sanskrit consonants.*

Like the Greeks, he drew his numbers from the Hindu script of alphabets; Sanskrit. To the 33 consonant in the Sanskrit script, he assigned a numerical value (see Fig.3). In Sanskrit, the first 25 letters are called Varga letters, the rest are called "Avarga" letters. Aryabhata assigned values, 1,2,3...25 to the 25 Varga letters. The avarga letters were assigned values 30,40,50,60,70,80,90 and 100 respectively. He also assigned numeric values to the nine Sanskrit vowels (see Fig.4).He then enunciated the rules to write down the numbers. If a consonant is in conjunction with a vowel, then the consonant value has to be multiplied by the place value of the vowel. Few examples of writing number in Aryabhata's numeric system are shown in Figure 4. One notes that Aryabhata could write down very large numbers. He was essentially an astronomer and as an astronomer, he needed large numbers. By his assignment of large numeric values to the Sanskrit vowels, he made sure that he can count very large numbers.

$$10^{16} \quad 10^{14} \quad 10^{12} \quad 10^{10} \quad 10^{8} \quad 10^{6} \quad 10^{4} \quad 10^{2} \quad 10^{0}$$
$$\text{औ} \quad \text{ओ} \quad \text{ऐ} \quad \text{ए} \quad \text{ऌ} \quad \text{ऋ} \quad \text{उ} \quad \text{इ} \quad \text{अ}$$

ग = 3
गु = 30000
गुण = 30015

*Figure 4: Numeric value assigned to nine sanskrit vowels. Example of writting numbers in Aryabhatta 's method is shown.*

Aryabhata wrote Aryabhatiya when he was only 23 years old. Thus, even around 500 CE, zero was not known to Indians. Zero was invented in India only around 628 CE. Zero has two uses, one as a number; symbol for absolute nothingness, and other as a placeholder. Certainly 2106 is not same as 216 or 2160. Indian mathematician Brahmagupta (598–c.670 CE) at the age of 30, wrote his magnum opus, "Brahmasphuṭasiddhanta", the first book that mentions zero as a number. As expected, he did not call it zero, rather he called it "sunya", a word still used in India, to indicate nothingness, vacuum. In the Brahmasphutasiddhanta, Brahmagupta defined zero as the result of subtracting a number from itself. He also gave the rules of mathematical operations with zero. As it was the custom of those days, the rules were given in a verse. If translated, they can be written as,

If a is a number, then,

(i)   $a+0=a$
(ii)  $a-0=a$
(iii) $a \times 0 = 0$
(iv)  $a/0 = a/0$
(v)   $0/0 = 0$

His rules of addition, subtraction and multiplication with zero are valid even today. However, his rules for division by zero were completely wrong. For the division of a number by zero, in many words, what Brahmagupta said amounts to a/0=a/0. Essentially, he was noncommittal about the division of a number by 0. Nowadays, we understand that the division by zero is not a mathematically defined operation, as there is no number which, multiplied by 0, gives a (assuming a≠0). However, sometimes it is useful to think a/0 as infinity. Even though Brahmagupta did conceive a symbol of nothingness, he could not conceive the notion of infinity. The concept of infinity is rather abstract, it is defined as something that is unlimited, endless, without bound. The concept can be casually understood as follows: imagine, you are counting numbers, 1,2,3,4..... , there is no end; you can go on counting one after another number. Indeed, whatever number you can count up to, I can always count a higher number just by adding 1. The number of natural numbers is then unlimited, endless. It is infinity. If you divide a number by zero again you get infinity. It can be understood as follow: divide a number continually with a diminishing number, the result will continue to increase, again without bound.

$$\frac{1}{.1}=10, \frac{1}{.001}=1000, \frac{1}{.00001}=10000, \frac{1}{.0000000001}=1000000000$$
$$\Rightarrow \frac{1}{0}=\infty \text{ (infinity)}$$

It may be noted here that infinity is not like any other number. One can add a number to infinity and it will be still be infinity. The difference between an ordinary number and infinity is best exemplified by giving the following example, which the great German mathematician David Hilbert (23 January 1862 – 14 February 1943) used to give to his students:

*In a hotel with a finite number of rooms, if fully occupied, manager has to refuse accommodation to a new guest. But not so in a hotel with infinite number of rooms. Say the rooms are numbered as; R1, R2, R3, R4, R5.... and it is fully occupied. How do the manger accommodate a new guest? He just shifts the existing guests to the next room;*

*Shift the guest at room R1 to room R2,*
*shift the guest at room R2 to room R3,*

*shift the guest at room R3 to room R4,*
*shift the guest at room R4 to room R5…*

*and so on. In the process, room R1 is released where the new guest can be accommodated.*

*If instead of single guest an infinite number of new guests arrive? Still no problem. The manager can do the following;*

*Shift the guest at room R1 to room R2,*
*Shift the guest at room R2 to room R4,*
*Shift the guest at room R3 to room R6,*
*Shift the guest at room R4 to room R8,*

*And so on. The process will release an infinite number of odd numbered rooms;*

True nature of infinity was understood only after George Cantor gave the set theoretic proof that there are infinite sets which cannot be put into one-to-one correspondence with the infinite set of natural numbers. Incidentally, the concept of infinity first appeared in ancient Greek. Greek philosopher Anaximander (ca. 610 – ca. 546 BCE), first took the giant leap in the mental thinking process when he conceived the "Apeiron", the infinite, unlimited or indefinite, along with eternal motion, as the originating cause of the world.

To come back to zero, Brahmagupta's rule division of zero by zero was also wrong. Since any number multiplied by zero is zero, the expression 0/0 also has no defined value. Indeed, as shown below, depending on the method of calculation, one can arrive at two equally valid answers,

$$\frac{.1}{.1}=1, \frac{.001}{.001}=1, \frac{.00001}{.00001}=1, \frac{.000000001}{.000000001}=1, \frac{.0000000000\mathit{1}}{.0000000000\mathit{1}}=1 \Rightarrow \frac{0}{0}=1,$$

$$\frac{0}{.1}=0, \frac{0}{.001}=0, \frac{0}{.00001}=0, \frac{0}{.000000001}=0, \frac{0}{.0000000000\mathit{1}}=0 \Rightarrow \frac{0}{0}=0,$$

Even though Brahmagupta, inventor of zero, could not give the correct rules for division of a number by zero, it is again in India, approximately after 500 years the correct rules emerged. The rules were given by the Indian mathematician Bhaskara II or Bhaskaracharya (1114-1185 CE). One of his mathematical treatises is Lilavati, named after his daughter. Bhaskaracharya understood the importance of zero, not just in positional notation, but also as a number. He gave the following eight rules concerning zero; for a non-zero number a,

$$(i)\ a \pm 0 = 0;\ (ii)\ 0^2 = 0;\ (iii)\ \sqrt{0} = 0;\ (iv)\ 0^3 = 0,$$
$$(v)\ \sqrt[3]{0} = 0;\ (vi)\ \frac{a}{0} = \infty;\ (vii)\ a \times 0 = 0;\ (viii)\ \frac{a \times 0}{0} = a$$

Bhaskaracharya understood infinity as a fraction who's denominator is 0. He did not call it infinity, rather he called it "Khahara" (खहर). In Lilavati, about the nature of Khahara or infinity, he wrote,

*"There is no change in Khahara (infinity) figure if something is added to or subtracted from the same."*

He then likened infinity with the almighty Vishnu,

*"It is like there is no change in infinite Vishnu (Almighty) due to dissolution or creation of abounding living beings."*

Incidentally, first, indubitable evidence of the use of zero in India, as early as 876 CE was found in a Vishnu temple, at Gwalior, a sprawling city of India. The temple has a stone tablet. For the temple, people of Gwalior gifted a garden measuring 187 by 280 hastas (One hasta approximately equals to 18 inches), to produce 50 garlands per day for worshipping. It was inscribed on the stone tablet. The inscribed 0 on the tablet is not much different from today's 0.

It may be mentioned here that it was also the Indian mathematician Brahmagupta who introduced the concept of negative numbers in mathematics. Greek mathematician Diophantus of Alexandria (200-ca.284 CE), who rightly is considered as the father of algebra, encountered negative numbers in the solution of linear equations. However, he disregarded them as unphysical, absurd. Europeans thought of numbers geometrically; a length, an area or a volume, all of which are real. Indeed in terms of length, area or volume, it is easy to understand;

2-1=1; or 2-2=0;

But what is 2-3=-1? How can you take a length 3 from a length of two? Indeed, even in 1803 the famous French mathematician Lazare Carnot was worried about the reality of negative numbers:

*"to really obtain an isolated negative quantity, it would be necessary to cut off an effective quantity from zero, to remove something of nothing: impossible operation. How thus to conceive an isolated negative quantity?"*

Unlike Europeans, Indians could think of numbers in an abstract manner. Brahmagupta, who first conceived the notion of absolute nothingness as a number zero, was also the first person in the history of mathematics to conceive negative numbers. He understood the negative numbers as a debt. He used two types of numbers, fortunes and debts, for positive and negative numbers respectively. He gave the following rules for working with negative numbers;

(1) A debt minus zero is a debt.
(2) A fortune minus zero is a fortune.
(3) Zero minus zero is a zero.
(4) A debt subtracted from zero is a fortune.
(5) A fortune subtracted from zero is a debt.
(6) The product of zero multiplied by a debt or fortune is zero.
(7) The product of zero multiplied by zero is zero.
(8) The product or quotient of two fortunes is one fortune.
(9) The product or quotient of two debts is one fortune.
(10) The product or quotient of a debt and a fortune is a debt.
(11) The product or quotient of a fortune and a debt is a debt.

From India, through the Arabian traders, Indian zero traveled to Arab. In the Golden age of Islam, during 800-1500 CE, Islamic scholars made spectacular progress in mathematics and astronomy. Caliph Harun-al-Rashid (reigned 786–809 CE) established an academy or intellectual centre at Bagdad called Bayt al-Ḥikmah or "House of Wisdom". House of Wisdom flourished during his son al-Ma'mun's regime, who himself was adept in various branches of science, medicine, philosophy, astronomy etc. Under his patronage, the pursuit of knowledge became a dominant feature of the Caliphate. Scholars were encouraged to translate medieval works and the vast Greek, and Indian literature, pertaining to philosophy, mathematics, natural science, and medicine were translated into the Arabic language. Renowned Arab mathematician, Al Khwarizmi (*ca.* 780- *ca* 850), studied in the House of Wisdom and later became its director. From his famous Book, Kitab al-jabr wa l-muqabala, the modern world got the term Algebra. Court of Bagdad received a gift of Brahmagupta's book and Al Khwarizmi overseered its translation. He understood the importance of zero as a number and place holder. He did not call is zero or sunya, rather he called if "sifr". In arabic, sifr means empty. Synthesizing the Arabian and Indian knowledge, in 825 CE, he published a book titled, *On the Calculation with Hindu Numerals.*

| Sanskrit | ० | ۹ | २ | ३ | ४ | ५ | ६ | ७ | ८ | ९ |
|---|---|---|---|---|---|---|---|---|---|---|
| Arabic | ٠ | ١ | ٢ | ٣ | ۴ | ۵ | ۶ | ٧ | ٨ | ٩ |
| English | 0 | 1 | 2 | 3 | 4 | 5 | 6 | 7 | 8 | 9 |

*Figure 5: The ten numerals of the decimal system in Sanskrit, Arabic and English.*

Zero entered Europe in the 12th century when Al Khwarizmi's book was translated into Latin under the title "al-Khwarizmi on the Numerals of the Indians". Italian mathematician Leonardo Fibonacci (ca. 1170 – ca. 1250 CE), also called Leonardo of Pisa, popularised the use of zero in Europe. Fibonacci is better known for introducing the integer sequence,

0, 1, 1, 2, 3, 5, 8, 13, 21, 34, 55, 89, 144, 233, 377.... .

The series is now known as Fibonacci series or sequence. The n-th term of the series is obtained by adding (n-1) and (n-2) terms. In 1202, Fibonacci wrote a book, Liber Abaci ( Book of Calculation). The book was instrumental in spreading use of zero in Europe. The book begins with,

*"The nine Indian figures are:*

*9 8 7 6 5 4 3 2 1.*

*With these nine figures and with sign 0 which the Arabs call zephir any number whatsoever is written, as is demonstrated below."*

The book was hugely popular in Europe and a second edition appeared in 1228. Slowly, Europeans accepted zero as a number and placeholder, and Italian "zephir" contracted to the present form "zero". How the ten decimal numbers evolved from India to Arab to Europe is schematically shown in Fig. 5. Traces of Indian origin can still be seen in few European numbers, European acceptance of the negative numbers took little longer. During the 15-16th century, Europeans reluctantly started to use negative numbers. Full and wholehearted acceptance came in the 17th century, only after the introduction of "Number line" by the English mathematician, John Wallis (1616-1703), who incidentally gave the present symbol of infinity ($\infty$), which represent an unending curve. In basic mathematics, a number line is a picture of an infinitely extended straight line on which every point corresponds to a real number. Any point in the line can be identified with zero. Then the points on the right of zero are positive numbers, the points on the left are negative.

For further reading few references are given below.